\documentclass[12pt]{amsart}

\usepackage{pstricks}
\usepackage[centertags]{amsmath}
\usepackage{amsfonts,amscd,latexsym}
\usepackage{pstcol}
\usepackage{euscript}
\usepackage[all]{xy}
\usepackage{showlabels}

\newcommand{\B}[1]{{\mathbb #1}}

\newtheorem{theorem}{Theorem}[section]

\theoremstyle{definition}

\newtheorem{example}[theorem]{Example}
\theoremstyle{remark}
\newtheorem{remark}[theorem]{Remark}
\newtheorem{question}[theorem]{Question}
\numberwithin{figure}{section}
\numberwithin{table}{section}

\newcommand{\al}{{\alpha}}

\newcommand{\om}{{\omega}}

\newcommand{\Mo}{(M,\omega )}

\newcommand\Symp{\operatorname{Symp}}

\newcommand{\NI}{{\noindent}}

\newpsobject{showgrid}{psgrid}{subgriddiv=2,griddots=10,gridlabels=5pt}

\begin{document}

\title{Fundamental group of $\Symp\Mo$ with no circle action}
\author{Jarek K\c edra}
\address{University of Szczecin, University of Aberdeen}
\curraddr{Mathematical Sciences\\
University of Aberdeen\\            
Meston Building\\
King's College\\
Aberdeen AB24 3UE\\
Scotland, UK.}
\email{kedra@maths.abdn.ac.uk}
\urladdr{http://www.maths.abdn.ac.uk/\~{}kedra}
\date{\today }
\thanks{Partially supported by the State Commitee for Scientific Reseach,
project number 1 PO3A 023 27.} 
\subjclass[2000]{Primary 53; Secondary 57}
\begin{abstract}
We show that $\pi_1(\Symp\Mo)$ can be nontrivial
for $M$ that does not admit any symplectic circle action.
\end{abstract}


\maketitle

\section{The question}

\begin{question}
Suppose that $\pi_1(\Symp\Mo)$ is nontrivial.
Is it true that some  nonzero element is represented
by a loop $S^1\to \Symp\Mo$ that is a homomorphism
(a circle action on $M$)?
\end{question}

\begin{remark}
If $G$ is a compact Lie group then 
any element of $\pi_1(G)$ is represented by a
loop that is a homomorphism.
Elements of $\pi_1(\Symp\Mo)$ which
are not represented by a circle action were described
in Anjos \cite{MR2003j:57052} as well as in McDuff
\cite{MR1929338}.
\end{remark}

\NI
{\bf Acknowledgements } I was asked the above question by Yael Karshon
during the conference in Stare Jab\l onki. The argument 
relies on a work of Lalonde and Pinsonnault \cite{MR2053755}.
I thank Rafa\l \, Walczak for comments.

\section{An answer}
We shall show that the answer to the above question is negative in
general. First we will give a little bit imprecise argument and
then show a concrete example. 

\subsection{A heuristic argument}
Philosophically speaking there are very few manifolds admitting a
circle action. They are, so to speak, exceptional. On the other
hand topology of groups of diffeomorphisms is rather complicated.
Hence one can expect nontrivial fundamental groups.

The argument consists of several steps:
\begin{enumerate}
\item
Take a closed simply connected symplectic manifold
$(K,\om_K)$.
Consider the evaluation fibration

$$\Symp(K,p)\to \Symp(K)\stackrel {ev}\to K,$$
where $\Symp(K,p) \subset \Symp(K)_0$ 
denote the isotropy subgroup of a point $p\in K$.

Observe that 
$ev_*:\pi_2(\Symp(K)) \to \pi_2(K)$
is trivial up to torsion. Indeed,
If $ev_*$ was nontorsion then
the corresponding map on rational cohomology would be nonzero,
say $ev^*(\al)\neq 0$ for $\al\in H^2(K,\B Q)$.
Then we have that $0=ev^*(\al^{n+1}) = ev^*(\al)^{n+1}$,
where $\dim K = 2n$.
But $\Symp(K)$ is a topological group so its rational
cohomology is free graded algebra. Thus every element of even
degree is of infinite order which contradicts the above calculation.
Thus we get that the rank of $\pi_1(\Symp(K,p))$ is not
smaller than the rank of $\pi_2(K)$. In particular it
is nonzero.
\item
The isotropy subgroup $\Symp(K,p)$ {\em should} be
weakly homotopy equivalent to the group of symplectomorphisms
of a one point blow-up of $(K,\om_K)$ in a very small ball.
This is proved for a range of 4-dimensional manifolds
by Lalonde and Pinsonnault in \cite{MR2053755}.
It is interesting to what extent it is true.
\item
The final step is to find a simply connected symplectic
closed manifold that its blow-up does not admit any symplectic
circle action.

\end{enumerate}

\subsection{Examples}

\begin{theorem}
Suppose that  $(K,\om_K)$ is neither rational
nor ruled surface up to blow-up.
Let $\Mo$ be a symplectic blow-up (in a small ball) of a 
closed simply connected K\"ahler surface $(K,\om_K)$. 
Then $\Mo$ admits no symplectic circle action and 
$\pi_1(\Symp\Mo)$ is nontrivial. 
\end{theorem}

\begin{proof}
We follow the scenario described in the previous section.
\begin{enumerate}
\item This step is obvious.
\item It is proved in \cite{K1} (Proof of Theorem 1.2) that $\pi_1(\Symp\Mo)$
surjects onto $\pi_1(\Symp(K,p))$. Hence $\pi_1(\Symp\Mo))$ is nontrivial.
The argument here relies on results of Lalonde and Pinsonnault
\cite{MR2053755}.
\item 
Suppose $\Mo$ admits a symplectic circle action. Since $\Mo$ is simply
connected then it has positive Euler characteristic. Hence the
action has a fixed point and i Hamiltonian. Thus, according to the classification of
Hamiltonian $S^1$-actions on  4-manifolds (see Theorem 13.19
in \cite{MR97b:58062}), we get that
$M$ has to be rational or ruled surface up to blow-up. This is excluded by the
hypothesis.

\end{enumerate}
\end{proof}

\begin{example}
To see a concrete example take $(K,\om_K)$ to be K3 surface with 
any symplectic form. Then its small blow-up $\Mo $ does not adtmit any
symplectic circle action and $\pi_1(\Symp\Mo)$ is nontrivial.
\end{example}

\bibliography{../../bib/bibliography}
\bibliographystyle{plain}

\end{document}